\documentclass[twoside, letterpaper, 12pt]{article}

\DeclareMathSizes{10}{10}{10}{10}
\usepackage{anysize}
\usepackage{fancyhdr}
\marginsize{2cm}{2cm}{0.5cm}{0.75cm}
\usepackage{setspace}
\usepackage[utf8]{inputenc} % allow utf-8 input
\usepackage[T1]{fontenc}    % use 8-bit T1 fonts
\usepackage{url}            % simple URL typesetting
\usepackage{booktabs}       % professional-quality tables
\usepackage{amsfonts}       % blackboard math symbols
\usepackage{nicefrac}       % compact symbols for 1/2, etc.
\usepackage{microtype}      % microtypography
\usepackage{amsthm}
\usepackage{caption}
\usepackage{subcaption}
\usepackage{graphics} % for pdf, bitmapped graphics files
\usepackage{epsfig} % for postscript graphics files
\usepackage{times} % assumes new font selection scheme installed
\usepackage{amsmath} % assumes amsmath package installed
\usepackage{amssymb}  % assumes amsmath package installed
\usepackage{makeidx}
\usepackage{float}
\usepackage{balance}
\usepackage{booktabs}
\usepackage{bbm}
\usepackage{mathrsfs}
\usepackage{enumerate}
\usepackage{multicol}
\usepackage{algorithm}
\usepackage{algpseudocode}

\usepackage{mathrsfs}  
\usepackage[normalem]{ulem}
\usepackage{xr} % enable referencing labels in other documents
\usepackage{multicol}
\usepackage{float}
\usepackage[nolist]{acronym}
\usepackage{xcolor}
\usepackage[colorlinks]{hyperref}       % hyperlinks

\hypersetup{
  colorlinks = true,
  allcolors = siaminlinkcolor,
  urlcolor = siamexlinkcolor,
}
\colorlet{siaminlinkcolor}{green!50!black}
\colorlet{siamexlinkcolor}{red!50!black}

\usepackage{setspace}
\usepackage{cleveref} %Cleveref must be defined after hyperref

\newcommand{\rd}{{\mathrm d}}

\newcommand{\vx}{{\bf x}}

\newcommand{\vPhi}{{\mbox{\boldmath$\Phi$}}}

\newcommand{\vmu}{{\mbox{\boldmath$\mu$}}}

\newcommand{\calH}{{\cal H}}

\newcommand{\calL}{{\cal L}}

\newcommand{\calU}{{\cal U}}

\newcommand{\argmin}{\operatornamewithlimits{argmin}}

\newcommand{\Rb}{\mathbb{R}}
\newcommand{\Eb}{\mathbb{E}}

\newcommand{\F}{\mathscr{F}}

\newcommand{\A}{\mathscr{A}}

\newcommand{\vm}{{\bf m}}

\newcommand{\vTheta}{{\mbox{\boldmath$\Theta$}}}

\makeatletter
\renewcommand{\maketag@@@}[1]{\hbox{\m@th\normalsize\normalfont#1}}%
\makeatother

\newtheorem{theorem}{Theorem}[section]

\newtheorem{prop}[theorem]{Proposition}

\catcode`\_=11\relax
    \newcommand\email[1]{\_email #1\q_nil}
    \def\_email#1@#2\q_nil{%
      \href{mailto:#1@#2}{{\emailfont #1\emailampersat #2}}
    }
    \newcommand\emailfont{}
    \newcommand\emailampersat{\small@}
    \catcode`\_=8\relax 
    
\makeatletter
\newcommand{\vast}{\bBigg@{3.5}}
\newcommand{\Vast}{\bBigg@{4}}
\newcommand{\vastt}{\bBigg@{4.5}}
\newcommand{\Vastt}{\bBigg@{5}}
\makeatother

\newcommand{\beginsupplement}{
        \setcounter{table}{0}
        \renewcommand{\thetable}{S\arabic{table}}
        \setcounter{figure}{0}
        \renewcommand{\thefigure}{S\arabic{figure}}
        \setcounter{section}{0}
        \renewcommand{\thesection}{S\arabic{section}}
        \setcounter{equation}{0}
        \renewcommand{\theequation}{S\arabic{equation}}
        \setcounter{algorithm}{0}
        \renewcommand{\thealgorithm}{S\arabic{algorithm}}
     }

\begin{acronym}
\acro{DNN}{Deep Neural Network}
\acro{ODE}{Ordinary Differential Equation}
\acro{SPDE}{Stochastic Partial Differential Equation}
\acro{FNN}{Feed-forward Neural Network}
\acro{CNN}{Convolutional Neural Network}
\acro{DP}{Dynamic Programming}
\acro{LSTM}{Long-Short Term Memory}
\acro{FC}{Fully Connected}
\acro{DDP}{Differential Dynamic Programming}
\acro{HJB}{Hamilton-Jacobi-Bellman}
\acro{PDE}{Partial Differential Equation}
\acro{PI}{Path Integral}
\acro{NN}{Neural Network}
\acro{SOC}{Stochastic Optimal Control}
\acro{RL}{Reinforcement Learning}
\acro{MPC}{Model Predictive Control}
\acro{IL}{Imitation Learning}
\acro{RNN}{Recurrent Neural Network}
\acro{DL}{Deep Learning}
\acro{SGD}{Stochastic Gradient Descent}
\acro{SDE}{Stochastic Differential Equation}
\acro{VRL}{Variational Reinforcement Learning}
\acro{IDVRL}{Infinite Dimensional Variational Reinforcement Learning}
\acro{1D}{1-dimensional}
\acro{2D}{2-dimensional}
\acro{ROM}{Reduced Order Model}
\acro{SDE}{Stochastic \ac{ODE}}
\end{acronym}

\begin{document}
	
	\title{
	\rule{\linewidth}{4pt}\vspace{0.3cm} \Large \textbf{
	  Variational Optimization Based Reinforcement Learning for Infinite Dimensional Stochastic Systems 
	}\\ \rule{\linewidth}{1.5pt}}
	\author{Ethan N. Evans$^{a,\dagger,}\thanks{Corresponding Author. Email: \email{eevans41@gatech.edu}}~$, Marcus A. Pereira$^{b,}\thanks{Authors contributed equally}~$, George I. Boutselis$^{a}$,\\ and Evangelos A. Theodorou$^{a,b}$\\ \vspace{-0.1cm}
	\small{$^a$Georgia Institute of Technology, Department of Aerospace Engineering} \\ \vspace{-0.2cm}
	\small{$^b$Georgia Institute of Technology, Institute of Robotics and Intelligent Machines} }
	
	\date{\small{This manuscript was compiled on \today}}
	
	\maketitle

\begin{abstract}
 
    Systems involving Partial Differential Equations (PDEs) have recently become more popular among the machine learning community. However prior methods usually treat infinite dimensional problems in finite dimensions with Reduced Order Models. This leads to committing to specific approximation schemes and subsequent derivation of control laws. Additionally, prior work does not consider spatio-temporal descriptions of noise that realistically represent the stochastic nature of physical systems. In this paper we suggest a new reinforcement learning framework that is mostly model-free for Stochastic PDEs with additive spacetime noise, based on variational optimization in infinite dimensions. In addition, our algorithm incorporates sparse representations that allow for efficient learning of feedback policies in high dimensions. We demonstrate the efficacy of the proposed approach with several simulated experiments on a variety of SPDEs.
  
\end{abstract}

% % Two or three meaningful keywords should be added here
% \keywords{Reinforcement Learning, Planning and Control}

\section{Introduction}
\label{sec:introduction}

Infinite dimensional stochastic  systems are typically systems that have spatio-temporal dynamics and are represented by \acp{SPDE}. Such systems appear in  areas of sciences and engineering such us fluid mechanics, plasma physics, partial observable stochastic control, continuum mechanics and many others.  Examples of  such systems are the stochastic Navier-Stokes equation which governs fluid flow and turbulence, the stochastic Nonlinear Schrödinger (NLS) equation, which has many realizations in Quantum Mechanics including the confinement of bosons in magnetic microtraps \cite{doi:10.1081/SAP-120017534}, the stochastic Nagumo equation which governs how voltage travels across a neuron in a brain \cite[Chapter 10]{lord_powell_shardlow_2014}, and the  stochastic Kuramoto-Sivashinsky (KS) equation which governs flame front propagation in combustion \cite{gomes2017controlling}. These fall into a category which covers a broad class of \acp{PDE} known as \textit{semi-linear} \acp{PDE}. A detailed exposition of certain examples in this category is given in  \cref{tab:semilinear_pdes}.

The theory of control of \ac{SPDE} systems was only introduced in the last few decades \cite{da1994stochastic, fabbri} and remains incomplete especially for the cases of stochastic boundary control. Numerical results and algorithms for distributed control of \acp{SPDE}  are limited and typically require some model reduction approach \cite{lou2009model,gomes2017controlling}.  In \cite{StochasticBurgers_1999}, the authors approach the control of the stochastic Burgers equation through the Hamilton-Jacobi-Bellman theory by applying the linear Feynman-Kac lemma; nevertheless, it lacks numerical results. In \cite{moura2013optimal}, the authors treat optimal control of linear deterministic \acp{PDE} by applying linear control theory, however this work is limited to linear \acp{PDE}. The book \cite{fabbri} gives a complete understanding of our ability so far, to apply optimal control theory to these systems.

\begin{table}[!t]
    \centering
    \begin{tabular}{ | l | p{2cm} | p{2.8cm} |p{3.5cm}|}
        \hline
        & \multicolumn{2}{c|}{\textbf{Operators}} & \\
        \hline  \textbf{Partial Differential Equation}
         & \textbf{Linear} $\A $  & \textbf{Non-linear} $F(t, X) $  & \textbf{State (or field)} \\ 
        \hline
        \hline
        \textcolor{orange}{Nagumo:} $u_t = {\color{blue}\epsilon u_{xx}} + {\color{violet}u(1-u)(u-\alpha)}$ & $u_{xx}$ & $u(1-u)(u-\alpha)$ & Voltage\\
        \hline
        \textcolor{orange}{Heat:} $u_t = {\color{blue}\epsilon u_{xx}}$ & $u_{xx}$ &  & Heat/temperature\\   \hline
        \textcolor{orange}{Burgers (viscous):} $u_t + {\color{violet}u u_x} = {\color{blue}\epsilon u_{xx}}$ & $u_{xx}$ & $uu_x$ & Velocity\\
        \hline
        \textcolor{orange}{Allen-Cahn:} $u_t = {\color{blue}\epsilon u_{xx} }+ {\color{violet}u - u^3}$ & $u_{xx}$ & $u - u^3$ & Phase of a material \\
        \hline
        \textcolor{orange}{NS:} $u_t = {\color{blue}\epsilon \Delta u} - \nabla p - {\color{violet}(u \cdot \nabla)u}$ & $\Delta u$ & $(u \cdot \nabla)u$ & Velocity\\        
        \hline
        \textcolor{orange}{NLS:} $iu_t + {\color{blue}\frac{1}{2}u_{xx}} + {\color{violet}|u|^2 u=0}$ & $u_{xx}$ & $|u|^2 u$ &  Wavefunction\\
        \hline 
        \textcolor{orange}{KdV:} $u_t + {\color{violet}6uu_x} + {\color{blue}u_{xxxx}}=0$ & $u_{xxxx}$ & $uu_x$ & Plasma wave \\
        \hline
        \textcolor{orange}{KS:} $u_t + {\color{violet}uu_x} + {\color{blue}u_{xx} + u_{xxxx}} = 0$ & $u_{xx} + u_{xxxx}$ & $uu_x$ & Flame front \\
        \hline
    \end{tabular} 
    \vspace{0.2cm}
    \caption{Examples of commonly known semi-linear PDEs in a \textit{fields representation} with $x$ representing spatial dimensions and $t$ representing time.}
    \label{tab:semilinear_pdes}
\end{table}

Contrasting with the work from the controls community are recent methods founded on machine learning techniques. These commonly treat deterministic \acp{PDE} as a finite set of \acp{ODE} through the use of \ac{ROM} type methods. Most recent work includes \cite{morton2018deep}, where the authors find reduced order Koopman-like local linear models and perform convex optimization for control using off-the-shelf solvers. However, this requires solving a least squares problem online which does not guarantee stabilizability of the resulting linear system. In \cite{rabault2019artificial} the authors successfully control a Navier-Stokes system with  reinforcement learning on policy networks in a deterministic, finite \ac{ODE} setting. Similarly, \cite{bieker2019deep} presents a Deep RNN framework with MPC to control a finite, deterministic \ac{ODE} representation (CFD solver) of a Navier-Stokes system.

\begin{figure}[!t]
\centering
  \includegraphics[width=0.5\linewidth]{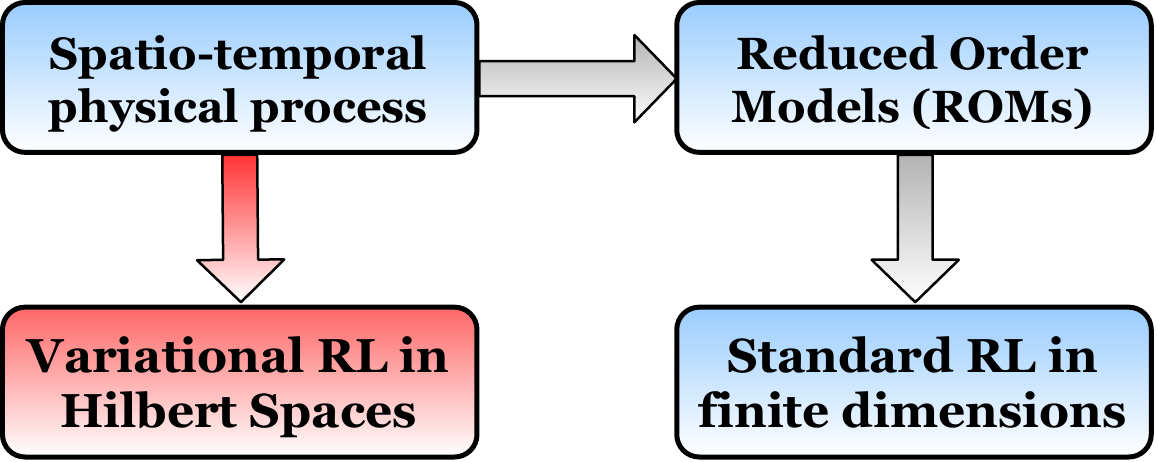}
  \caption{Our proposed approach versus traditional approaches.}
\label{fig:approach_diagram}
\end{figure}

In contrast to recent work which first require developing deterministic \acp{ROM} and then using standard approaches from \ac{RL} or \ac{MPC}, we treat the \ac{SPDE} system directly in Hilbert spaces, and derive a variational optimization framework for episodic reinforcement of policy networks as highlighted in red in \cref{fig:approach_diagram}.

We take inspiration from a general principle stemming from statistical physics and thermodynamics that has been shown to have applicability in stochastic optimal control \cite{TheodorouCDC2012}:
\begin{equation}\label{eq:Free_Energy_Relative_Entropy}
\text{Free Energy} \leq \text{Work} - \text{Temperature} \times \text{Entropy} 
\end{equation}
Optimization of this relation from a measure theoretic perspective gives rise to the well known Gibbs measure which is used in variational inference problems \cite{wainwright2008graphical}. This perspective enables us to seek a middle ground between recent results in \ac{DL} and traditional stochastic optimal control: We approach \acp{PDE} with infinite dimensional stochastic calculus, yet apply highly successful \ac{DL} techniques. We develop a new method fusing together variational optimization, episodic reinforcement learning, and measure theoretic stochastic calculus in infinite dimensions.

\vspace{1em}
Our specific contributions are as follows:
\begin{itemize}
    \item We present a reinforcement learning framework in Hilbert spaces based on variational optimization and importance sampling for \acp{SPDE}. The resulting algorithm  incorporates explicit feedback of the entire \ac{SPDE} and allows for arbitrary non-linear polcies such as \acp{FNN}, \acp{CNN} and \acp{RNN}
    \item We introduce a technique to handle numerical integration of policy networks over the spatial domain  which we call $SparseForwardPass$ for \ac{FNN} and \ac{CNN} policies, enabling scalability to 2D and 3D problems.
    \item Since the algorithm is derived in infinite-dimensional space, any choice of numerical approximation scheme such as finite difference, spectral Galerkin or finite-element can be used to approximate trajectory samples. In addition, due to performing optimization in infinite dimensional space, the derivation is valid for the stochastic versions of all PDEs included in \cref{tab:semilinear_pdes} and therefore is general. \end{itemize}

\section{Problem Formulation}
\label{sec:formulation}

%===============================================================================
% Zakai: $\rd u = (\frac{1}{2}\partial_{xx} (a u) - \partial_x(bu))\rd t -\partial_x (\gamma u) \rd W$

This work proposes control of a large class of infinite-dimensional systems described by \acp{SPDE} that are of \textit{semi-linear} form. There are other ways to express such systems, however here we take the approach of expressing the system as evolving on time-indexed separable Hilbert spaces in order to leverage several mathematical tools developed in such spaces. Consider the general semi-linear controlled \ac{SPDE} given by 
\begin{align} \label{eq:SPDEs_Control}
\rd X &= \big( \A X   + F(t, X)  \big) \rd t  + G(t, X)\big(\Phi(t,X,\vx; \vTheta^{(k)})\rd t+ \frac{1}{\sqrt{\rho}} \rd  W(t)\big), 
\end{align}
where $X(t) \in \calH$ is the state of the system which evolves on the Hilbert space $\calH$, the linear and nonlinear operators $\A: \calH \rightarrow \calH$ and $F(t,X): \Rb \times \calH \rightarrow \calH$  (resp.) are uncontrolled drift terms, $\Phi(t,X,\vx; \vTheta^{(k)}):\Rb \times \calH \times \Rb^3\rightarrow \calH$ is the nonlinear control policy paramterized by $\vTheta^{(k)}$ at the $k^{th}$ iteration, $\rd W(t):\Rb \rightarrow \calH$  is a Cylindrical spatio-temporal noise process (i.e. space-time white noise), and $G(t,X)$ is nonlinearity that affects both the Cylindrical noise and the control. It is used to incorporate the effects of actuation on either the field (distributed) or at the boundaries. Referring back to \cref{tab:semilinear_pdes}, the generality of the \textit{Hilbert spaces formulation} becomes clear as any semi-linear \ac{PDE} can be handled by appropriately choosing $\A$ and $F$. For a more complete introduction, including some mild but necessary assumptions and clear definitions of the Cylindrical process, see the supplementary material and the references therein.

Define the uncontrolled and controlled probability measures associated with \cref{eq:SPDEs_Control} as $\mathcal{L}$ and $\tilde{\calL}$, respectively. These measures roughly describe the probabilistic evolution of the system, with the probability density function as a finite dimensional analog. In this case, \cref{eq:Free_Energy_Relative_Entropy} takes the form
\begin{align} \label{eq:Legendre}
  & - \frac{1}{\rho}   \log \Eb_{\calL} \bigg[ \exp( -\rho {J} )  \bigg]  = \min_{\calU(\cdot,\cdot)} \bigg[    \Eb_{\tilde{\calL}}\left({J} \right)  + \frac{1}{\rho} D_{KL} ( \tilde{\calL}\hspace{0.05cm} \big|\big|  \calL )  \bigg],
\end{align}
where $J=J(X)$ can be viewed as an arbitrary state cost function. The associated ``Work" and ``Entropy" terms that minimize this expression describe a minimum ``energy"\footnote{The term energy here is used loosely to describe the landscape for work and entropy} measure. Sampling from this measure would simultaneously minimize state cost and the $KL$-divergence between the controlled and uncontrolled distributions, which in this case is roughly interpreted as control effort. The measure that optimizes \cref{eq:Legendre} is the so-called Gibbs measure
\begin{equation}\label{eq:Gibbs}
\rd \calL^{*} = \frac{\exp( - \rho J) \rd \calL}{\Eb_\calL \big[\exp( - \rho J) \big] }.
\end{equation}

% \vspace{-0.2cm}
While it is not known how to sample directly from \cref{eq:Gibbs}, the goal of variational optimization methods is to incrementally reduce the distance (defined in the Kullback–Leibler divergence sense) between the controlled distribution $\tilde{\calL}$ and the optimal measure \cref{eq:Gibbs}. We formulate our variational minimization problem as 
\begin{align} \label{eq:theta}
    \vTheta^{*} &=  \argmin_{\vTheta}  D_{KL}(\calL^{*}|| \tilde{\calL}) \nonumber \\
    &= \argmin_{\vTheta} \bigg[\int \log \Big(\frac{\rd \calL}{\rd \tilde{\calL}} \Big)  \frac{\rd  \calL^*}{\rd \calL} \frac{\rd  \calL}{\rd \tilde{\calL}} \rd \tilde{\calL} \bigg] = \argmin_{\vTheta} \;L
  \end{align}
A more detailed derivation can be found in the supplementary. Finally, we introduce a version of Girsanov's Change of Measure theorem (found in supplementary) between the uncontrolled and controlled processes, resulting in the so-called Radon-Nikodym derivative given as
\small\begin{align}
\label{eq:radon_i_nonlinear_feedback}
\begin{split}
\frac{ \rd \calL}{ \rd \tilde{\calL} } = \exp\bigg\lbrace-\sqrt{\rho}  \int_{0}^{T} \Big\langle \vPhi(t,X,\vx ; \vTheta^{(k)}), \rd W(t) \Big\rangle   -\rho\frac{1}{2}\int_{0}^{T} \Big\langle \vPhi(t,X,\vx ; \vTheta^{(k)}), \vPhi(t,X,\vx ; \vTheta^{(k)}) \Big\rangle \rd t\bigg\rbrace
\end{split}
\end{align} \normalsize
Plugging in \cref{eq:Gibbs} and  \cref{eq:radon_i_nonlinear_feedback} (for importance sampling), the loss-function $L$ becomes
\begin{align}\label{eq:loss_function}
    L =  \Eb_{\tilde{\calL}}  \vast[\underbrace{ \frac{\exp( - \rho \tilde{J})}{\Eb_{\tilde{\calL}} \big[ \exp( - \rho \tilde{J}) \big]}}_{Importance Weight} \bigg( &-\sqrt{\rho} \int_{0}^{T} \underbrace{ \Big\langle \vPhi(t,X,\vx ; \vTheta^{(k)}), \rd W(t) \Big\rangle}_{Noise Inner Product
} \nonumber \\
    & - \frac{1}{2} \rho \int_{0}^{T} \underbrace{\Big\langle \vPhi(t,X,\vx ; \vTheta^{(k)}), \vPhi(t,X,\vx ; \vTheta^{(k)}) \Big\rangle}_{Policy Inner Product} \rd t \bigg) \vast],
\end{align}
where $\tilde{J}$ is defined by
\begin{equation}\label{eq:importance_sampled_cost}
    \tilde{J} = \underbrace{ \vphantom{\Big\langle \vPhi(t,X,\vx ; \vTheta^{(k)}), \rd W(t) \Big\rangle}J}_{State Cost} + \frac{1}{\sqrt{\rho}}\int_{0}^{T} \underbrace{\Big\langle \vPhi(t,X,\vx ; \vTheta^{(k)}), \rd W(t) \Big\rangle}_{Noise Inner Product}  +\frac{1}{2}\int_{0}^{T} \underbrace{\Big\langle \vPhi(t,X,\vx ; \vTheta^{(k)}), \vPhi(t,X,\vx ; \vTheta^{(k)}) \Big\rangle}_{Policy Inner Product} \rd t
\end{equation}
The intermediate steps that lead to the above final forms of \cref{eq:theta} and \cref{eq:loss_function} can be found in the supplementary material.
The loss-function $L$ exponentiates the cost of the system trajectories, evaluated by $\tilde{J}$, to produce a weighted average of the mixed control-noise term and the quadratic control term. We minimize this loss via \ac{SGD}. The resulting Variational RL with learn rate $\gamma$ is an incremental update of the form 
\begin{equation}\label{eq:parameter_update}
\begin{split}
    \vTheta^{(k+1)} =  \vTheta^{(k)} - \gamma \nabla_{\vTheta}L.
\end{split}
\end{equation}

In contrast to prior work that use variational optimization to approximate optimal probability measures such as \cite{Williams2016AggressiveDW}, our proposed approach uses an arbitrary non-linear feedback policy as compared to time-varying policy of step-functions. Such policies also lead to different parameter update-rules requiring inversion of a jacobian. \ac{SGD}-based minimization also lends itself to the use of well-known backprop-based algorithms such as ADA-Grad \cite{duchi2011adaptive} and  ADAM \cite{kingma2014adam}.

Although the state may be described by an infinite-dimensional vector in a Hilbert space, many physical realizations of actuation are defined on finite subspaces. The above derivation keeps $\vPhi$ as mapping into the Hilbert space, insinuating that the actuation may be distributed everywhere and infinite-dimensional. However, the goal of this work is to ultimately use finite-action policy networks to control \cref{eq:SPDEs_Control}. As such, we refine $\vPhi$ as
\begin{equation}
    \vPhi(t,X,\vx;\vTheta^{(k)}) = \vm(\vx)^\top \varphi(X ; \vTheta^{(k)} )
\end{equation}
where $\varphi(X; \vTheta^{(k)}): \calH \rightarrow \Rb^N$ is a finite policy network with $N$ control outputs representing $N$ distributed (or boundary) actuators. The function $\vm(\vx):D \rightarrow \Rb^N \times \calH$ represents the effect of the finite actuation on the infinite-dimensional field, where $D$ is the domain of the finite spatial region. Some examples of $\vm(\vx)$ are Gaussians-like exponential functions with mean centered at the actuator location (for distributed control) and indicator functions (for boundary control).

%===============================================================================

\section{Algorithm and Network Architecture}\label{sec:algorithm}

% Above we present the derivation of the variation reinforcement learning method in Hilbert spaces. The authors conjecture that beyond mathematical elegance, there is an algorithmic benefit to performing optimization in Hilbert spaces before finally discretizing the field for implementation on a computer. In the field of computational PDEs, there are several different ways one can spatially discretize a PDE. However, the \ac{IDVRL} algorithm is agnostic to the spatial discretization. This is in contrast to reduced order methods which discretize spatially and temporally before deriving a suitable method. Each derived method in the reduced order approach depends heavily on the style of discretization.

The above derivation provides a mathematical framework for updating the weights of a policy network in a \ac{RL} setting. In order to implement it as an algorithm, data must be generated either from a physics-based or data-based model, or from interactions with a real system. Notice that since the only term from the dynamics to appear in  \cref{eq:loss_function,eq:importance_sampled_cost} is the Cylindrical noise term $\rd W$, there is no need to have an explicit \ac{SPDE} model. As a result, any black-box methods that incorporate spatio-temporal stochasticity can be used to generate sample trajectories of the system.

\begin{figure}[!ht]
\centering
  \includegraphics[width=0.75\linewidth]{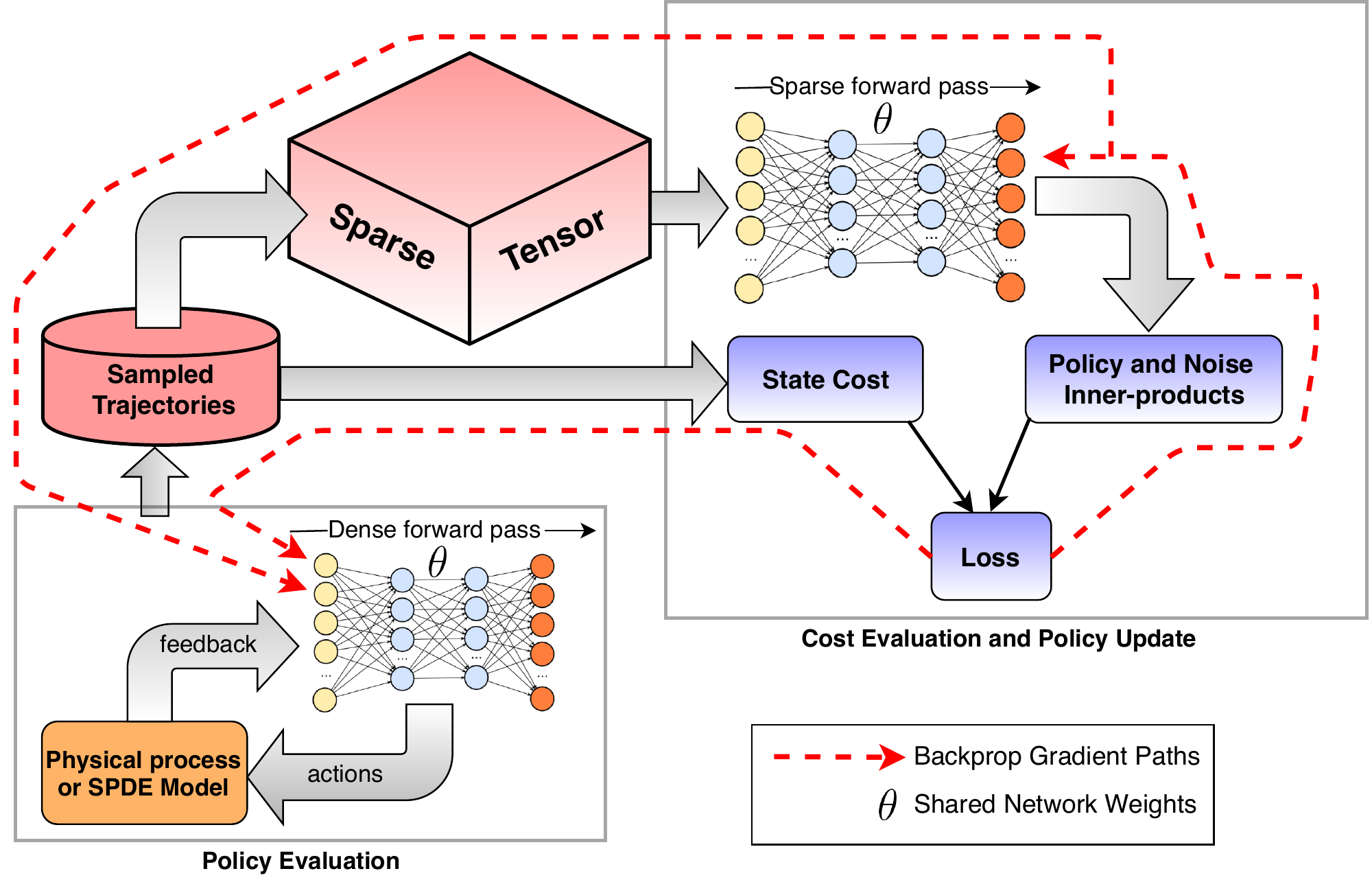}
  \caption{Block diagram of computational graph for Variational Reinforcement Learning of Infinite-dimensional systems.} %State trajectories are generated based on either an {SPDE} model or  black-box model with the policy network providing actuation. A sparse tensor is generated and a sparse forward pass enables computation of $PolicyInnerProduct$ and $NoiseInnerProduct$. Sampled trajectories are also used to compute state cost, resulting in a computation of the loss. Finally gradient descent produces backpropagation through various paths resulting in an update to the network parameters.}
\label{fig:architecture_diagram}
\end{figure}

The above derivation introduces a unique problem for our proposed reinforcement learning framework that has not been addressed in prior work. Each inner product in Hilbert space in \cref{eq:loss_function,eq:importance_sampled_cost} represents a spatial integration over a finite region $D$. To the knowledge of the authors, integration over a policy network has not been attempted to date. However in this work, we integrate spatially over the input to the network. Consider the inner product indicated as \textit{PolicyInnerProduct}. The representation of this inner product as a spatial integration takes the form 
\small \begin{align} \label{eq:inner_product_computation}
    \int_{0}^{T} \Big\langle \vPhi(X,\vx ; \vTheta^{(k)}), \vPhi(X,\vx ; \vTheta^{(k)}) \Big\rangle \rd t &= \int_{0}^{T} \iint_D \varphi(X(t,x,y) ; \vTheta^{(k)} )^\top M(x,y) \varphi(X(t,x,y) ; \vTheta^{(k)} ) \rd x \rd y   \rd t \nonumber \\
    &= \int_{0}^{T} \sum_{j=1}^{\infty} \varphi(X(e_j) ; \vTheta^{(k)} )^\top M(e_j) \varphi(X(e_j) ; \vTheta^{(k)} ) \rd t,
\end{align}\normalsize
where $D \subseteq \Rb^2$ is the problem domain, $\lbrace e_j \in \calH : j = 0,1,2,\dots \rbrace$ forms an orthonormal basis over $\calH$, and $M(\vx) = \vm(\vx)\vm(\vx)^\top$. After discretization on a 2D grid, the basis becomes a finite set $\lbrace e_j \in \Rb^{J^2} : j = 0,1,2,\dots \rbrace$, where each element is a one-hot vector. Thus, evaluating the spatial integral is reduced to summing up forward passes through the policy network with each pixel considered individually.  

Spatially integrating over the policy network is a memory intensive task, where the storage becomes ($J^2, J, J$) for each sample over the time horizon. However, given that the basis elements of each $(J,J)$ ``image" have only one activated ``pixel", the resulting tensor is tremendously sparse. As such, each layer's activation can be computed with a sparse matrix multiplication, resulting in what we call a $SparseForwardPass$ method that is not memory intensive for relatively large 2D problems. This can be applied to both \acp{FNN} and \acp{CNN}. For \acp{CNN}, activation can be achieved by matrix multiplication with a Toeplitz matrix constructed from the filter coefficients \cite{chellapilla2006high}.

A summary of our architecture is depicted in \cref{fig:architecture_diagram}. A policy network with initialized weights is passed through a model or physical realization of the system to produce state trajectories, which are used to compute a state cost as well as a sparse tensor that is used to compute the inner products in \cref{eq:loss_function,eq:importance_sampled_cost} in a memory and time-efficient manner. Finally the loss is computed and passed to a gradient-based optimization algorithm. This approach is independent of specific policy network architecture used, which can often be problem dependent. In this work we used two different networks: a \ac{FNN} for \ac{1D} \ac{SPDE} and a \ac{CNN} for \ac{2D} \ac{SPDE}.

The resulting \ac{IDVRL} algorithm is shown in  \cref{Algorithm1}, wherein subscript implies an element of the corresponding vector. The input terms are time horizon ($T$), number of iterations ($K$), number of rollouts ($R$), initial state ($X_0$), number of actuators ($N$), noise variance ($\rho$), time discretization ($\Delta t$), actuator locations ($\mu$), actuator variance ($\sigma_\mu$, for distributed control cases), and initial network parameters ($\vTheta^{(0)}$). We note that the function $GradientOptimize(L, \vTheta^{(k)})$ represents the update from \cref{eq:parameter_update}. As mentioned above, this is handled by any variant of \ac{SGD}, which performs backpropagation through the network. The computational graph of the proposed algorithm has multiple backprop paths, as shown by the dotted red line in \cref{fig:architecture_diagram}. For more information on $SampleNoise()$, refer to \cite[Chapter 10]{lord_powell_shardlow_2014}. 

\begin{algorithm}[!t]
 \caption{\acl{IDVRL}}
 \begin{algorithmic}[1]
 \State \textbf{Function:} \textit{$\vTheta^* =$ \textbf{OptimizePolicyNetwork}($T,K,R,X_0,N,\rho,\Delta t,\mu, \sigma_{\mu},\vTheta^{(0)}$)}
 \State Compute $\vm(\vx), M(\vx)$ $\forall$ $\vx \in D$
 \For {$k=1 \; \text{to} \; K$}
 \For {$r=1 \; \text{to}\; R$}
 \For{$t=1 \;\text{to}\; T$}
     \State $\rd W_t \gets SampleNoise()$
     \State $X_t \gets Propagate(X_{t-1},\vTheta^{(k)}, \rd W_t)$ via \cref{eq:SPDEs_Control}
     \State $J_r \gets J_r + StateCost(X_t)$ 
     \State $S_t \gets SparseForwardPass(\vTheta^{(k)},X_t)$
     \State $N_t \gets NoiseInnerProduct\big(S_t, \rd W_t,  \vm(\vx)\big)$
     \State $P_t \gets PolicyInnerProduct\big(S_t, M(\vx)\big)$
 \EndFor
 \State $P,\; N \gets Sum(P_t), \; Sum(N_t)$
  \State $\tilde{J}_r \gets \tilde{J}(P,N,J_r)$
%  \State $N \gets Sum(N_t)$
 \EndFor
 \State $W \gets ImportanceWeight(\tilde{J})$
 \State $L \gets ComputeLoss(P,N,W)$ via \cref{eq:importance_sampled_cost}
 \State $\vTheta^{(k+1)} \gets GradientOptimize(L,\vTheta^{(k)})$
 \EndFor
 \end{algorithmic}
 \label{Algorithm1}
\end{algorithm}

%===============================================================================

\section{Simulation Results and Discussion}
\label{sec:simulation}

We now briefly discuss our simulations and results for both distributed and boundary control tasks on 3 different \acp{SPDE}. All tasks were reaching tasks, wherein the policy has to achieve a desired profile in certain parts of the spatial domain. The computational graph for all simulations was implemented in Tensorflow \cite{tensorflow} to leverage GPU parallelization for training as well as sparse linear algebra operations for $SparseForwardPass$. We used a time-discretization of $\Delta t=0.01\,s$ and $\Delta t=0.02\,s$ for \ac{1D} and \ac{2D} problems respectively. The data for training the policies was generated by simulating the \acp{SPDE} using centered finite-difference approximation for the spatial derivatives on a \ac{1D} or \ac{2D} grid and a semi-implicit Euler scheme for discretization of the time derivatives. For detailed explanation on these schemes, we refer the reader to \cite[Chapters 3 and 10]{lord_powell_shardlow_2014}. For 1D simulations, we used an Alienware laptop with an Intel Core i9-8950HK CPU @ 2.9GHz$\times$12, $32$ GB RAM and a NVIDIA GeForce GTX 1080 graphics card. On average, Tensorflow-GPU required around $16$ minutes of training time for $1000$ iterations. For the 2D simulation, we used Tensorflow-CPU, due to insufficiency of VRAM, which required around $12$ hours of training time for $1000$ iterations.  

\begin{figure}[!t]

\begin{multicols}{3}
    \begin{subfigure}[h!]{0.33\textwidth}
    \includegraphics[width=\textwidth]{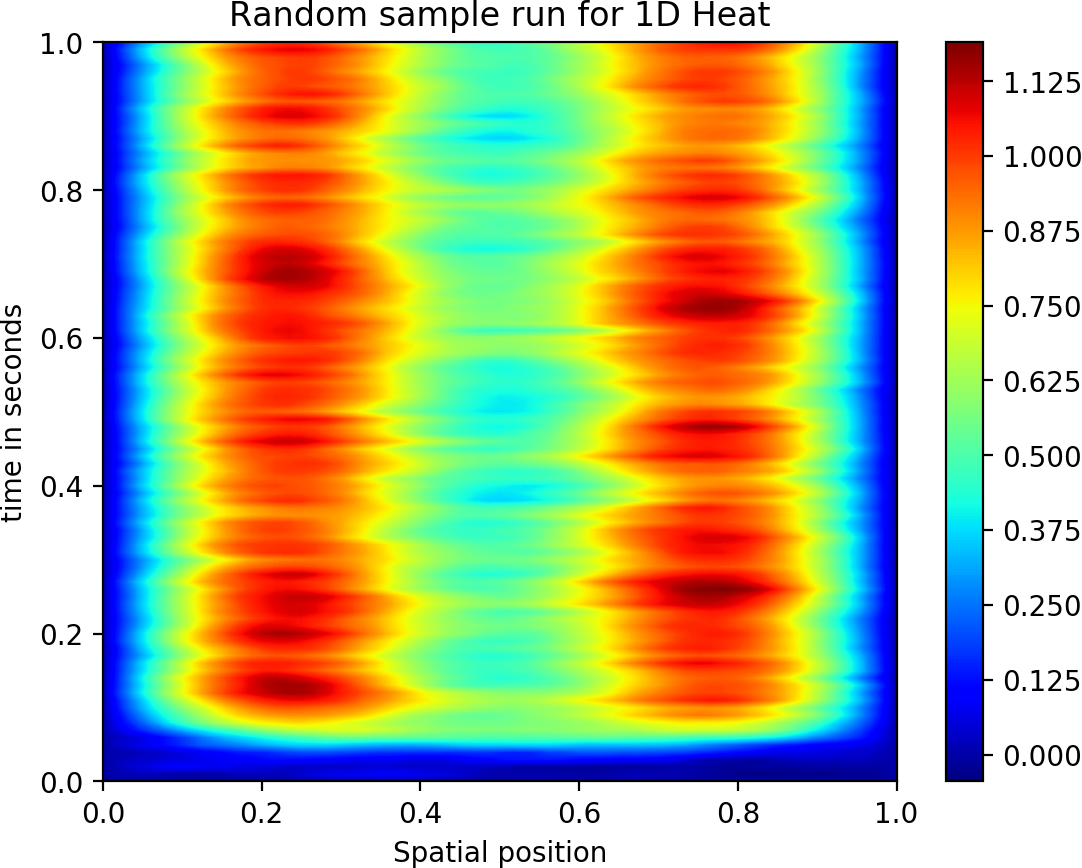}
    \label{fig:Heat1D_contour}
    \vspace{-0.5cm}
    \caption{}
    \end{subfigure}
    
    \begin{subfigure}[h!]{0.33\textwidth} \includegraphics[width=\textwidth]{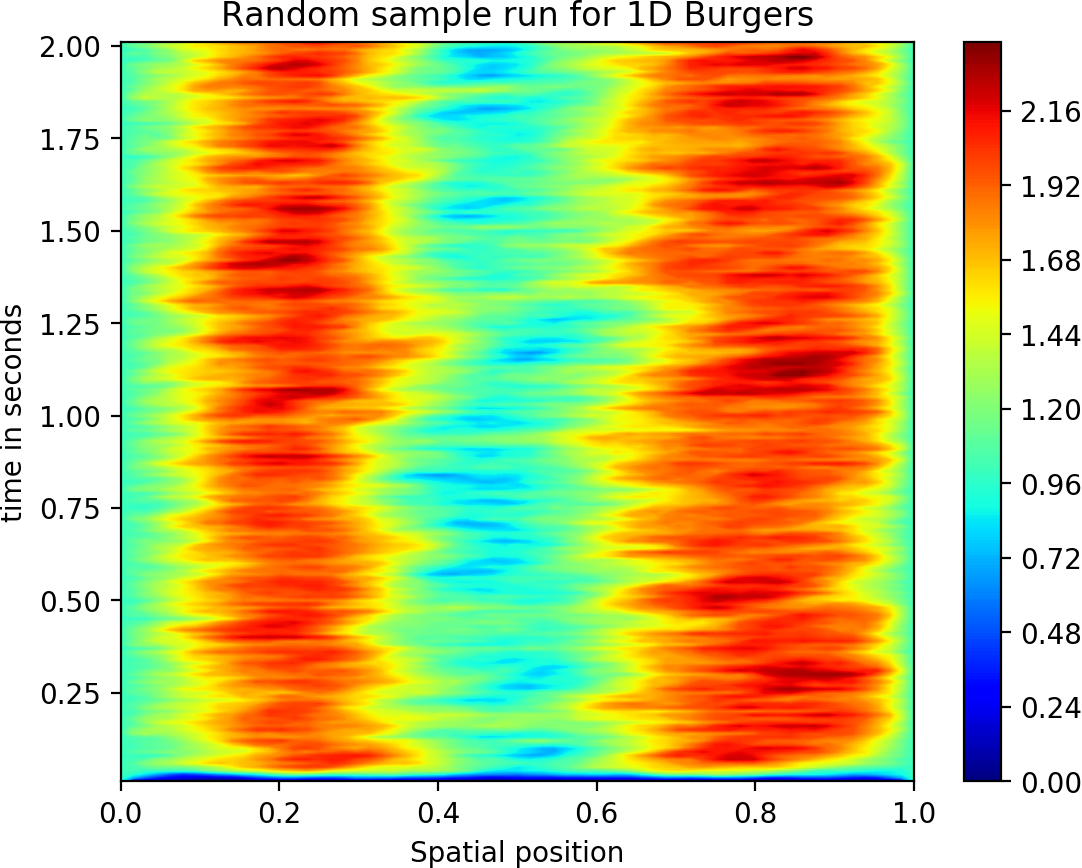}
    \vspace{-0.5cm}
    \caption{}
    \label{fig:Burgers1D_contour}
    \end{subfigure}
    
    \begin{subfigure}[h!]{0.33\textwidth} \includegraphics[width=\textwidth]{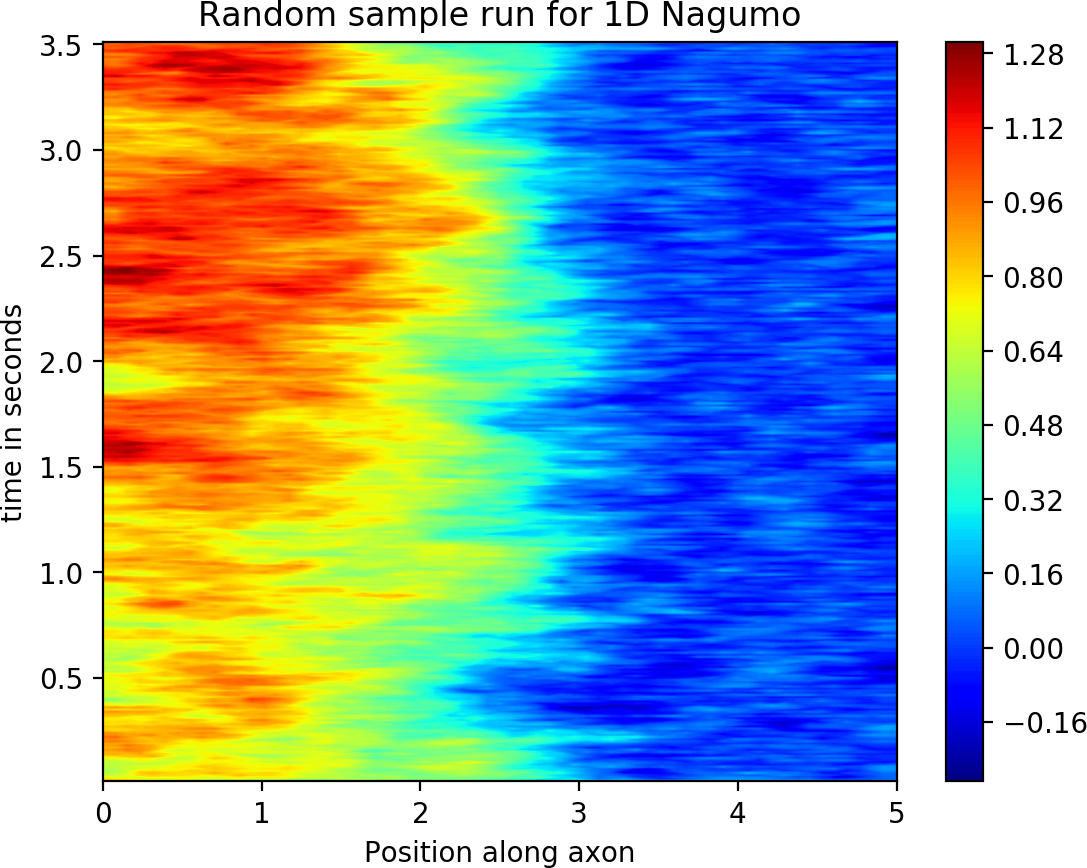}
    \vspace{-0.5cm}
    \caption{}
    \label{fig:Nagumo1D_contour}
    \end{subfigure}            

\end{multicols}

\vspace{-0.75cm}\begin{multicols}{3}
    \begin{subfigure}[h!]{0.365\textwidth}    \includegraphics[width=\textwidth]{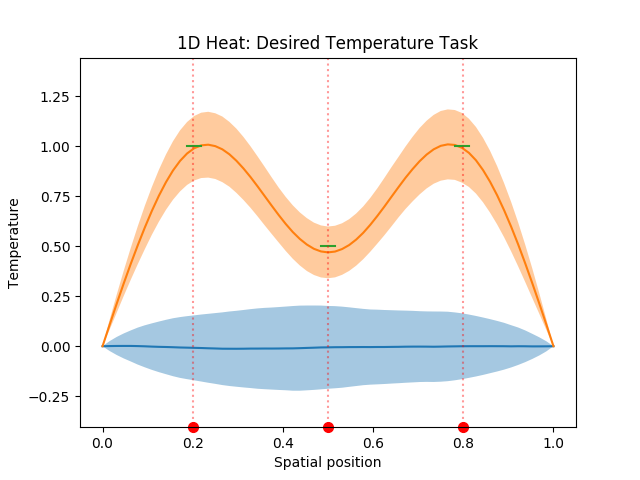}
    \vspace{-0.5cm}
    \caption{}
    \label{fig:Heat1D_traj}
    \end{subfigure}
    
    \hspace{-0.15cm}\begin{subfigure}[h!]{0.365\textwidth}
   \includegraphics[width=\textwidth]{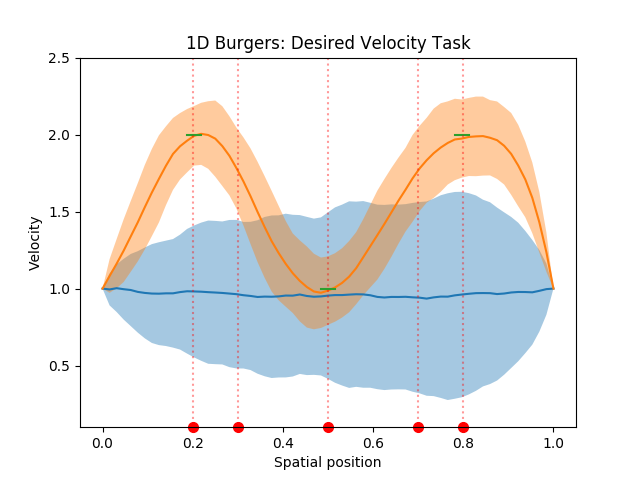}
    \vspace{-0.5cm}
    \caption{}
    \label{fig:Burgers1D_traj}
    \end{subfigure}
     
    \hspace{-0.3cm}\begin{subfigure}[h!]{0.365\textwidth}
    \includegraphics[width=\textwidth]{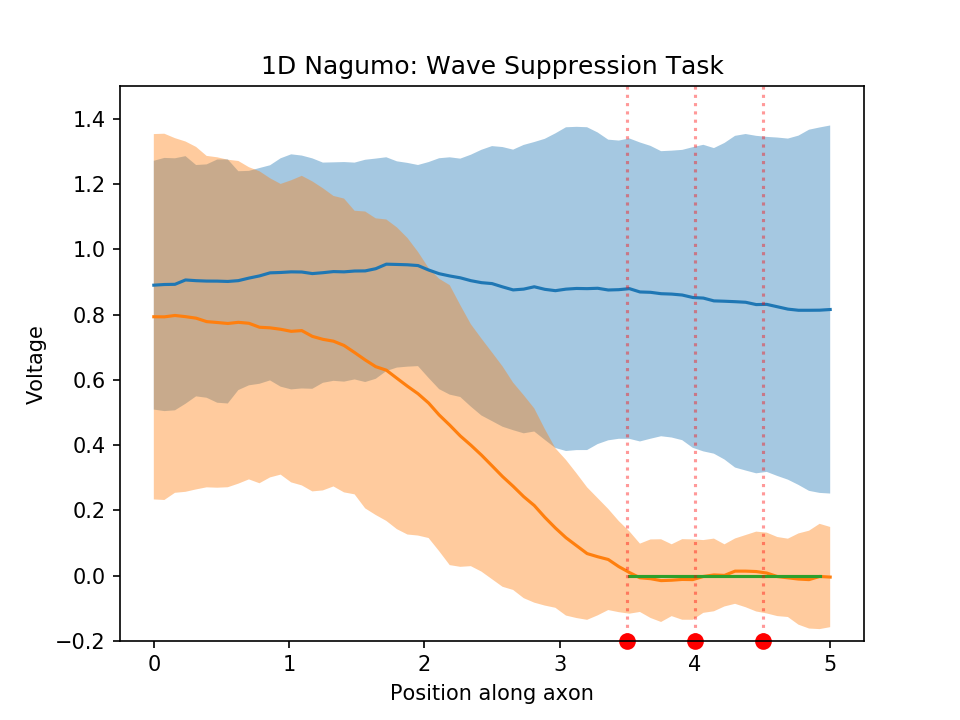}
    \vspace{-0.5cm}
    \caption{}
    \label{fig:Nagumo1D_traj}
    \end{subfigure}        
\end{multicols}

\centering
\begin{multicols}{2}
        \hspace{1.75cm}\begin{subfigure}[h!]{0.33\textwidth} \includegraphics[width=\textwidth]{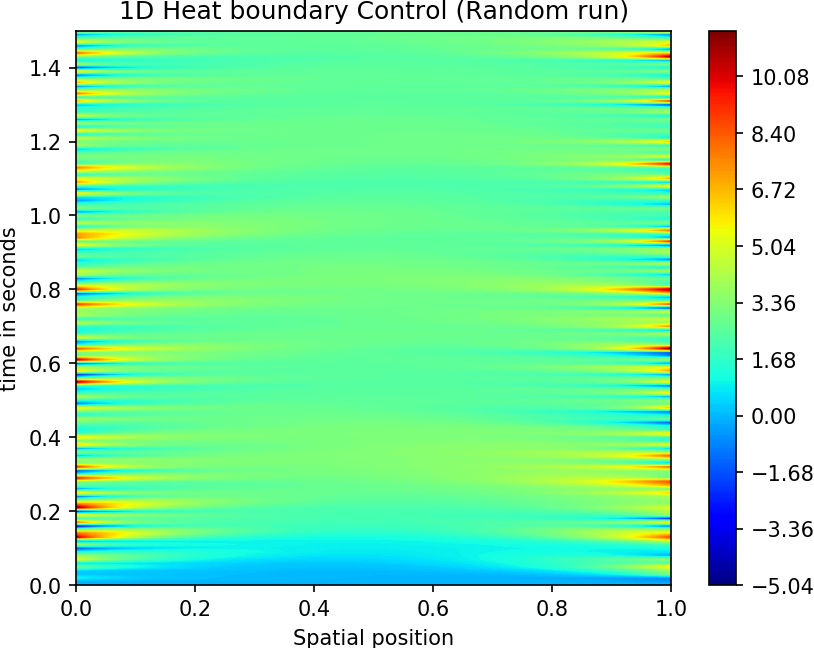}
        \vspace{-0.5cm}
        \caption{}
        \label{fig:Heat1D_boundary_contour}
        \end{subfigure}            
        
        \hspace{-1.5cm}\begin{subfigure}[h!]{0.365\textwidth} \vspace{-0.3cm}\includegraphics[width=\textwidth]{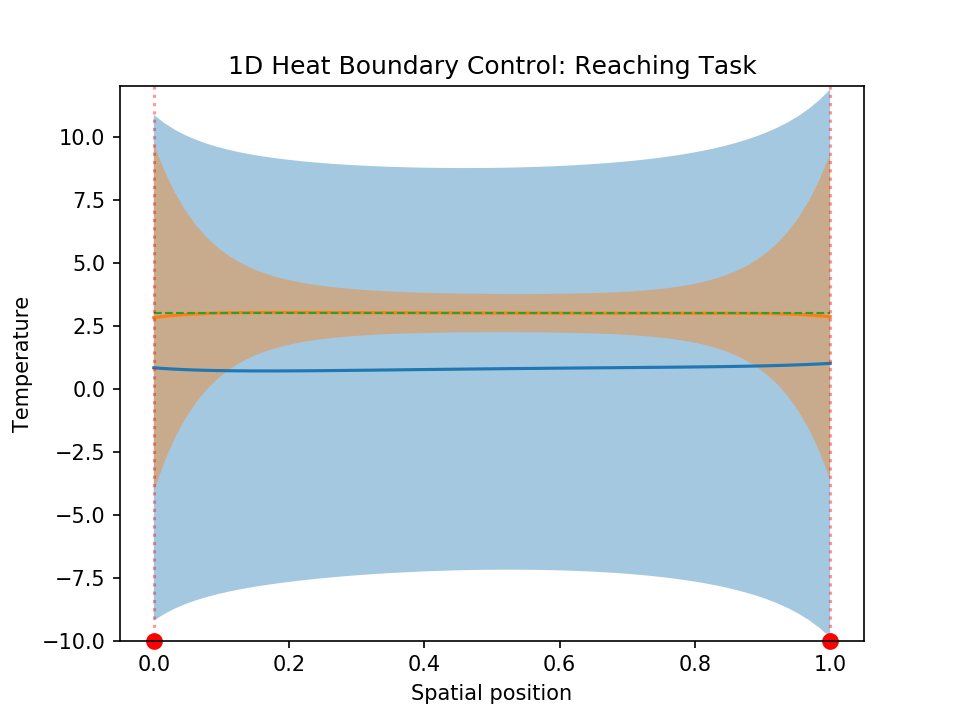}
        \caption{}
        \label{fig:Heat1D_boundary_traj}
        \end{subfigure}     
\end{multicols}

\vspace{-0.25cm}\caption{\textbf{Control of \ac{1D} \acp{SPDE}.} (a), (d), (g), (h) correspond to the Heat \ac{SPDE}, (b), (e) to Burgers \ac{SPDE}, and (c), (f) to Nagumo \ac{SPDE}. In (d), (e), (f), (h)  we use blue to represent \textit{mean uncontrolled profiles}, orange to represent \textit{mean controlled profiles} generated by the trained policy network, green to represent \textit{desired values} along certain parts of the spatial domain and red to represent \textit{locations of actuator centers}. The means and variances represent statistics at the final simulation time over 200 trials. (a), (b), (c), (g) depict a randomly picked trial run to emphasize the presence of spatio-temporal stochasticity. (a-f) depict results for distributed control of \acp{SPDE} and (g-h) depict results for boundary control of a \ac{SPDE}} 
\label{fig:1D_pdes}
\end{figure}

Figure \ref{fig:1D_pdes} (a) and (d) depict the results of the \ac{IDVRL} algorithm on a task of controlling the \ac{1D} heat \ac{SPDE} with homogeneous Dirichlet boundary conditions. The goal of the task is to raise and maintain the temperature to $T=1$ at regions around $x=0.2$ and $x=0.8$, and $T=0.5$ at a region around $x=0.5$. Figure \ref{fig:1D_pdes}a) shows the temperature contours of a single realization of the completed task and \cref{fig:1D_pdes}d) shows the mean controlled and uncontrolled trajectories at the final time with a 2-$\sigma$ variance shaded in the corresponding color. The boundary conditions fixed the endpoints to a temperature of $T=0$, as shown. 

Figure \ref{fig:1D_pdes}, (b) and (e) depict the results of the \ac{IDVRL} algorithm on the task of controlling the \ac{1D} Burgers \ac{SPDE} with non-homogeneous Dirichlet boundary conditions. In this task the goal is to reach a desired velocity in the medium at given locations. This is challenging given the nonlinear advection behavior of the system in addition to the pure diffusion behavior shown in the \ac{1D} heat \ac{SPDE} task. The advection-diffusion creates an apparent rightwards wave-front that must be accounted for by the policy network in order to achieve the task. Given the increased difficulty of the problem, we added actuators, as indicated by vertical red dotted lines. Despite the added actuators, the task remains severely under-actuated.

Figure \ref{fig:1D_pdes}, (c) and (f) depict the \ac{IDVRL} algorithm on the task of controlling the \ac{1D} Nagumo \ac{SPDE} with homogeneous Neumann boundary conditions. As noted earlier, the Nagumo SPDE represents voltage travelling across the axon of a neuron in the brain. The goal of this task is to suppress the voltage from travelling across the axon. Voltage near $1.0$ indicates the voltage has travelled across, and in this suppression task, we seek to keep the voltage at the right end of the axon at $V=0$. As shown in \cref{tab:semilinear_pdes}, the Nagumo SPDE has a 3rd order nonlinearity. For this task, we supplied the system with only three actuators near the right end, where voltage must be suppressed.

\begin{figure*}[t]
        \begin{multicols}{4}
            \begin{subfigure}[h!]{0.235\textwidth}
            \includegraphics[width=\textwidth]{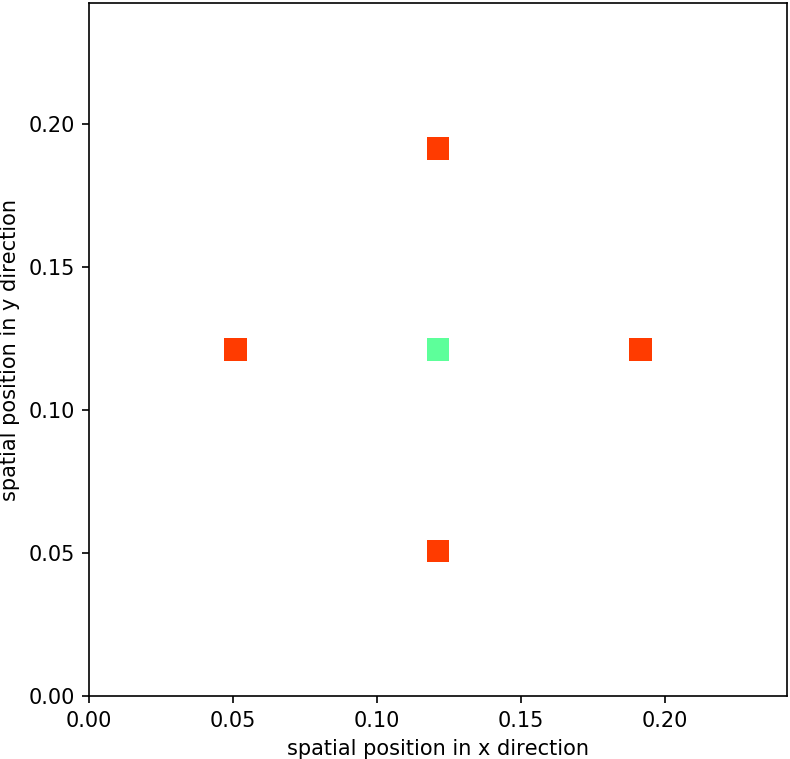}
            \vspace{-0.5cm}
            \caption{}
            \label{fig:2Dheat_actuators}
            \end{subfigure}
            
            \begin{subfigure}[h!]{0.225\textwidth} \includegraphics[width=\textwidth]{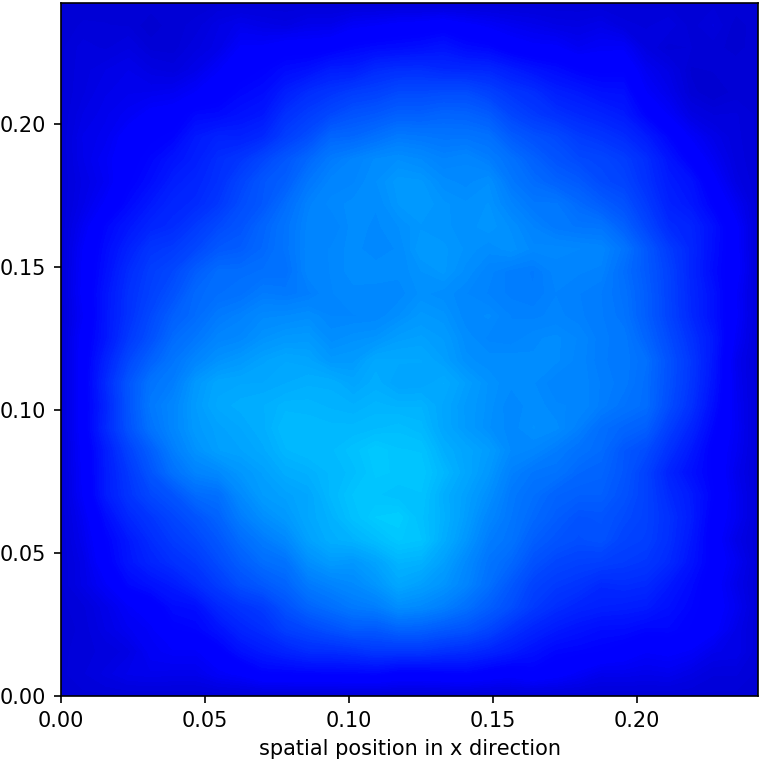}
            \vspace{-0.5cm}
            \caption{}
            \label{fig:2Dheat_IC}
            \end{subfigure}
            
            \begin{subfigure}[h!]{0.225\textwidth} \includegraphics[width=\textwidth]{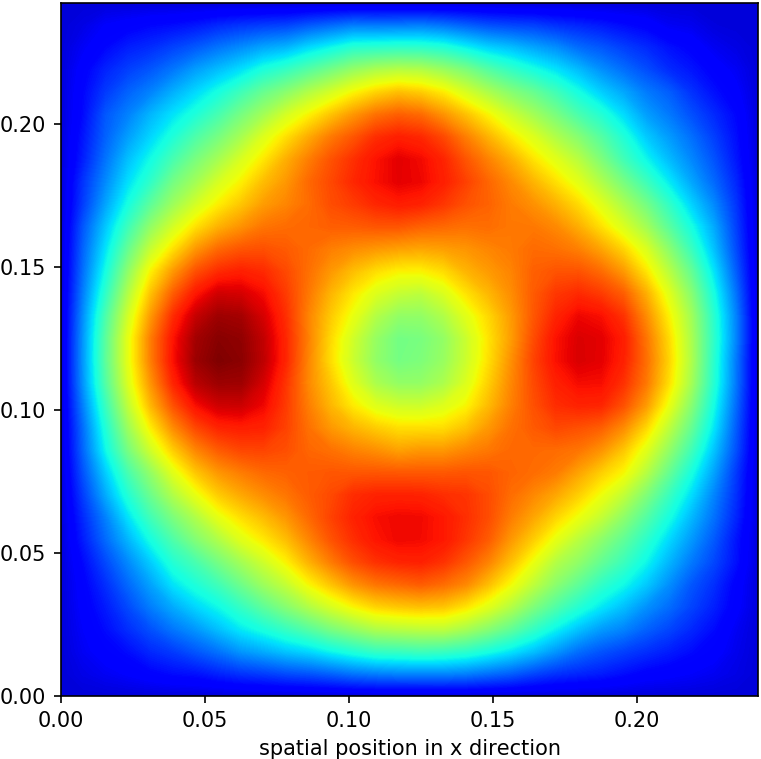}
            \vspace{-0.5cm}
            \caption{}
            \label{fig:2Dheat_middle}
            \end{subfigure}            

            \begin{subfigure}[h!]{0.245\textwidth} \includegraphics[width=\textwidth]{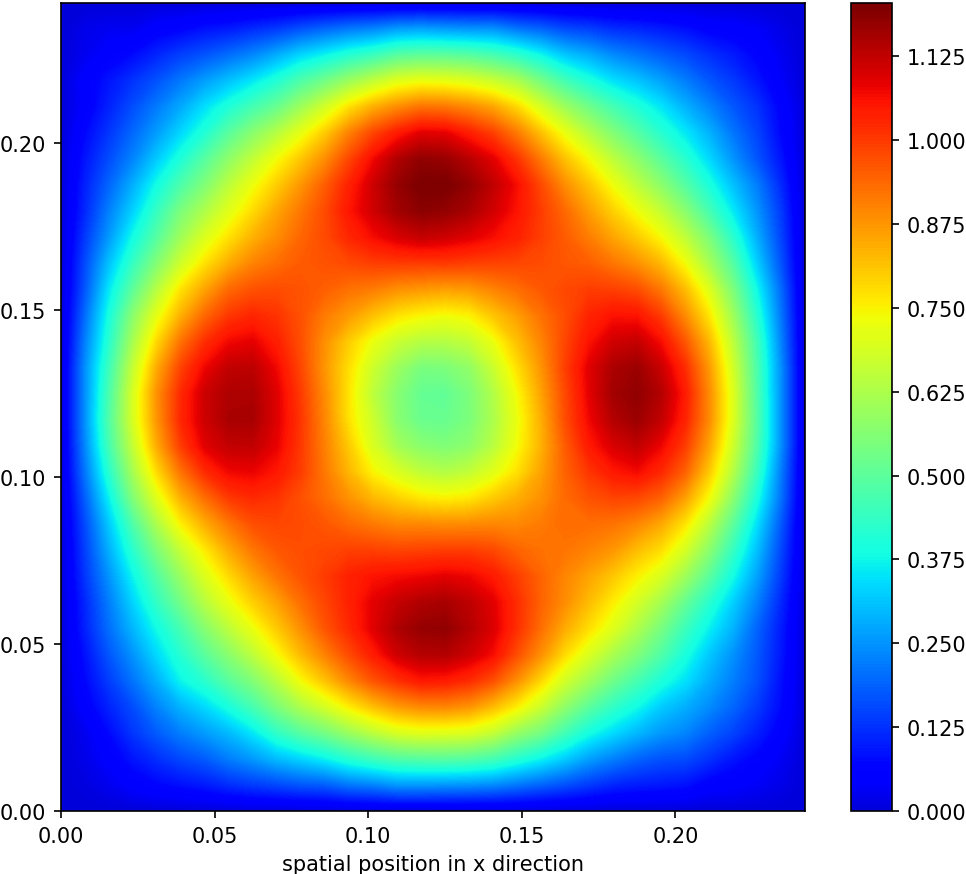}
            \vspace{-0.5cm}
            \caption{}
            \label{fig:2Dheat_final}
            \end{subfigure} 
            
        \end{multicols}
        \vspace{-0.25cm}\caption{\textbf{Control of the 2D Heat \ac{SPDE}}. The plot on the extreme \textit{left} shows the desired profile patches for the reaching task. Next is the profile at the start of simulation, followed by profile half-way through and on the extreme \textit{right} is the profile at end of simulation time. This is a randomly picked instance from a batch of test trials. The color-bar depicts the range of temperatures in the simulation. The locations of the 5 actuator-centers are at the centers of the desired patches.}
        \label{fig:2DHeat}
\end{figure*}

For the next task, we scaled the \ac{IDVRL} algorithm to two-dimensional problems. With this task we attempt to control the \ac{2D} Heat \ac{SPDE} with homogeneous Dirichlet boundary conditions with a \ac{CNN} policy network. The goal of this task it to raise the temperature in five regions. The desired temperature at the four outer regions is $T=1$ and the desired temperature at the center region is $T=0.5$. Figure \ref{fig:2DHeat} depicts a single realization of the controlled task under a significant amount of noise with five actuators.

In contrast to the previous tasks where actuators are distributed in the field, \cref{fig:Heat1D_boundary_traj}) depicts a \textit{boundary} control task, where the actuator controls the boundary condition. The Radon-Nikodym derivative exists for the case of boundary control of semi-linear \acp{SPDE} with boundary noise \cite{da1994stochastic}, and we demonstrate that our method similarly extends to this case. The task here is similar to the first case, where the policy network is tasked with reaching a desired value of $T=3$.
   
We invite the interested reader to refer to our supplementary material for specific details on each of our simulations such as cost functions, hyper-parameter values, neural network parameters and videos comparing controlled and uncontrolled \acp{SPDE}.

Throughout our simulated experiments, especially for distributed control tasks, we found that the algorithm is not sensitive to the majority of our parameters. We noted that a useful heuristic in applying the algorithm to new problems without having to tune the parameters was to ensure that the starting loss function was not very close to zero (i.e. 1e-10). Despite a large variance of noise that we typically applied to our systems ($\rho = 10$), the optimization algorithm was able to converge in under 1000 iterations for \ac{1D} problems and under 2000 iterations for \ac{2D} problems. 

On the whole, even though injecting higher variance noise into the system inherently makes the control task much more challenging, high variance noise is useful in our algorithm for exploration over rollouts at each iteration. As such, there is an inverse relationship for a given convergence behavior between variance in the noise and number of rollouts.

% \begin{figure}[t]
%     \begin{multicols}{2}
%         \begin{subfigure}[h!]{0.45\textwidth} \includegraphics[width=\textwidth]{Figures/heat1d/heat_1d_boundary_contour.png}
%         \vspace{-0.5cm}
%         \caption{}
%         \end{subfigure}            

%         \hspace{-0.25cm}\begin{subfigure}[h!]{0.525\textwidth} \vspace{-0.3cm}\includegraphics[width=\textwidth]{Figures/heat1d/heat_1d_boundary.png}
%         \vspace{-0.65cm}
%         \caption{}
%         \end{subfigure}     
%     \end{multicols}
%     \vspace{-0.25cm}\caption{\textbf{1D Heat Boundary Control.}  (a) \textit{mean uncontrolled profiles}, orange to represent \textit{mean controlled profiles} generated by the trained policy network, green to represent \textit{desired values} along certain parts of the spatial domain and red to represent \textit{locations of actuator centers}. The means and variances represent statistics at the final simulation time over 200 trials. (a), (b), (c) depict a randomly picked trial run to emphasi }
%     \label{fig:heat_1d_boundary}
% \end{figure}

There are also several interesting behaviors that the \ac{IDVRL} algorithm demonstrates. First, we noticed that often times throughout optimization, the loss would decrease as desired, but state cost would temporarily increase, before decreasing more dramatically after some number of iterations. This indicates that there may not be a strictly proportional relationship between loss function and state cost. Indeed a lower state cost implies that the task is being accomplished, yet a trend of decreasing loss function indicated that when there was a temporary increase in state cost, the \ac{IDVRL} algorithm may have been pushing the network parameters out of a local minimum towards better task performance in later iterations. These trends indicate that the \ac{IDVRL} algorithm may perform well on experiments outside the ones considered in this paper.

%===============================================================================

\section{Conclusion and Future Directions}
\label{sec:conclusion}

This work presents a variational reinforcement learning algorithm for the distributed and boundary control of infinite dimensional stochastic systems. The optimization method was derived in Hilbert spaces, thereby avoiding the need to depend on specific descretization schemes to realize the algorithm. The resulting algorithm requires only an actuation model and therefore is mostly model-free. The algorithm was demonstrated on five simulated experiments including \ac{1D} and \ac{2D} with both distributed and boundary type actuation.

In future work the authors will investigate provable convergence properties for \ac{IDVRL} based on \cite{zhou2014gradient}, and will implement the algorithm on some \acp{SPDE} described in this paper such as the Stochastic Navier-Stokes equation using state-of-the art CFD solvers. The authors also plan to investigate second-order SPDEs such as the Euler-Bernoulli equation which has been used to investigate the dynamics of tentacle-like soft continuum robots \cite{jin2015lyapunov}.

%===============================================================================

% Bibliography
%\section*{References}

\bibliographystyle{ieeetran}
\balance
\bibliography{References}

%===============================================================================

\newpage
\section*{Supplementary Information}
\beginsupplement

\section{Q-Wiener and Cylindrical noise} \label{app:q_wiener}
The Q-Wiener noise defines a generalization of the standard Brownian motion for fields. The cylindrical noise which is used in our paper is a specific case of a Q-Wiener process. Below we give an important property of the aforementioned noise profiles which makes a connection to the standard Wiener noise. Formal definitions can be bound in \cite[Chapter 4]{da1992stochastic}.

\begin{prop} Let $ \lbrace e_{i}\rbrace_{i=1}^{\infty} $  be a complete orthonormal system for the Hilbert Space $U$  such that $ Q e_{i}= \lambda_{i} e_{i}   $. Here, $ \lambda_{i} $ is the eigenvalue of $ Q$ that corresponds to  eigenvector $ e_{i} $.  Then,  a   $Q$-Wiener process $ W(t) \in U$ has the following expansion:
 \begin{equation}\label{supeq:Q_Wiener}
    W(t) = \sum_{j=1}^{\infty} \sqrt{\lambda_{j}} \beta_{j}(t) e_{j},
 \end{equation}
\noindent where  $ \beta_{j}(t)  $  are real valued Brownian motions that are mutually independent on $ (\Omega, \F, \mathbb{P}) $.
\end{prop}

\noindent We note that when $\lambda_{j}=1$ $\forall j\in \mathbb{N}$, $W(t)$ corresponds to a cylindrical Wiener process (space-time white noise). In this case, the series in \eqref{supeq:Q_Wiener} converges in another Hilbert space $U_{1}\supset U$, when the inclusion $\iota:U\rightarrow U_{1}$ is Hilbert-Schmidt. For more details see \cite{da1992stochastic}.

\section{Girsanov's theorem and the Randon-Nikodym derivative}
\begin{theorem}[Girsanov] \label{girs} Let $\Omega$ be a sample space with a $\sigma$-algebra $\mathcal{F}$. Consider the following $H$-valued stochastic processes:
\begin{align}
\rd X&=\big(\A X+F(t, X)\big) \rd t+G(t, X)\rd W(t), \label{X}\\
\rd\tilde{X}&=\big(\A \tilde{X}+F(t, \tilde{X})\big)\rd t+\tilde{B}(t, \tilde{X})\rd t+G(t, \tilde{X})\rd W(t),\label{eq:eqsss}
\end{align}
where $X(0)=\tilde{X}(0)=x$ and $W\in U$ is a cylindrical Wiener  process with respect to measure $\mathbb{P}$. Moreover, let $\Gamma$ be a set of continuous-time, infinite-dimensional trajectories in the time interval $[0,T]$. Now the {\it probability law} of $X$ will be defined as $\mathcal{L}(\Gamma):=\mathbb{P}(\omega\in\Omega|X(\cdot,\omega)\in\Gamma)$ . Similarly, the law of $\tilde{X}$ is defined as $\tilde{\calL}(\Gamma):=\mathbb{P}(\omega\in\Omega|\tilde{X}(\cdot,\omega)\in\Gamma)$. Then
\begin{equation} \label{eq:Lg}
\begin{split}
\tilde{\calL}(\Gamma)=  \textstyle\mathbb{E}_{\mathbb{P}}\big[\exp\big(\int_{0}^{T}\langle\psi(s), dW(s)\rangle_{U}-\frac{1}{2}\int_{0}^{T}||\psi(s)||_{U}^{2}ds\big)|X(\cdot)\in\Gamma\big],
\end{split}
\end{equation}
where we have defined $\psi(t):=G^{-1}\big(t, X(t)\big)\tilde{B}\big(t, X(t)\big)\in U_{0}$ and assumed
%\begin{equation}
 %   \label{psi_assum1}
  $  \mathbb{E}_{\mathbb{P}}\big[e^{\frac{1}{2}\int_{0}^{T}||\psi(t)||^2\mathrm dt}\big]<+\infty$.
%\end{equation}% Here, we write for brevity $\tilde{\calL}(\omega)\equiv\tilde{\calL}(\tilde{X}(\cdot,\omega)\in\Gamma)$.
\end{theorem}
The proof of Girsanov's theorem can be found in \cite{DaPrato1999}. It follows that the {\it Randon-Nikodym} (RN) derivative between measures $\cal L(\cdot)$ and $\tilde{\calL}(\cdot)$ of the different dynamical systems defined in \eqref{eq:eqsss}, is given by
\begin{equation} \label{eq:RN}
\begin{split}
\frac{\mathrm d\tilde{\calL}}{\mathrm d{\calL}} = \exp\big(\int_{0}^{T}\langle\psi(s), dW(s)\rangle_{U}-\frac{1}{2}\int_{0}^{T}||\psi(s)||_{U}^{2}ds\big),
\end{split}
\end{equation}
In the main paper, we use the RN derivative for the case where $\cal L(\cdot)$, $\tilde{\calL}(\cdot)$ correspond to uncontrolled and controlled trajectories, respectively, with $\psi(\cdot)$ being properly defined.

%Let $H$, $U$ be separable Hilbert spaces and let $(\Omega, \F, \Pb) $ be a probability space with filtration $\F_{t},$ $t\in [0, T] $. The uncontrolled and controlled dynamics are (respectively) given by
%\begin{align}
%\rd X &= \big(\A X   + F(t, X)\big)  \rd t +   \frac{1}{\sqrt{\rho}} G(t, X) \rd W(t), \label{supeq:SPDEs_NoControl}\\
%\rd \tilde{X} &= \big( \A \tilde{X}   + F(t, \tilde{X})  \big) \rd t  + G(t, \tilde X)\big(\calU(t, \tilde X)\rd t+ \frac{1}{\sqrt{\rho}} \rd  W(t)\big), \label{supeq:SPDEs_Control}
%\end{align}
%with initial conditions $X(0) =\tilde{X}(0)= \xi$. Here, $W\in U$ is a cylindrical Wiener process on $(\Omega, \mathcal{F}, \mathbb{P})$.

\section{Derivation of Variational Minimization and Loss Function}

This section explains the steps to arrive at equations \cref{eq:theta,eq:loss_function,eq:importance_sampled_cost} from the main paper. 

\begin{align} \label{supeq:theta}
    \vTheta^{*} &=  \argmin_{\vTheta}  D_{KL}(\calL^{*}|| \tilde{\calL}) \nonumber \\
    &= \argmin_{\vTheta} \bigg[   \int  \log\Big( \frac{ \rd \calL^{*}}{\rd \tilde{\calL}} \Big)  \rd  \calL^{*} \bigg] \nonumber \\
    &= \argmin_{\vTheta} \bigg[  \int \log \Big( \frac{ \rd \calL^*}{\rd \calL} \frac{\rd \calL}{\rd \tilde{\calL}} \Big)  \rd  \calL^* \bigg] \nonumber \\
    &=  \int \log \Big(\frac{ \rd \calL^*}{\rd \calL}\Big) \rd \calL^* + \argmin_{\vTheta} \bigg[ \int \log \Big( \frac{\rd \calL}{\rd \tilde{\calL}} \Big) \rd  \calL^* \bigg] \nonumber \\
    &= \argmin_{\vTheta} \bigg[\int \log \Big(\frac{\rd \calL}{\rd \tilde{\calL}} \Big)  \frac{\rd  \calL^*}{\rd \calL} \frac{\rd  \calL}{\rd \tilde{\calL}} \rd \tilde{\calL} \bigg] = \argmin_{\vTheta} \;L
  \end{align}
Now,
\begin{align*}\label{suppeq:loss_function}
    L &= \Eb_{\tilde{\calL}}  \Bigg[\log \Big(\frac{\rd \calL}{\rd \tilde{\calL}} \Big)  \frac{\rd  \calL^*}{\rd \calL} \frac{\rd  \calL}{\rd \tilde{\calL}} \rd \tilde{\calL} \Bigg]
\end{align*}
Substituting \cref{eq:radon_i_nonlinear_feedback}, the log goes away because of the exponential,
\begin{align*}
    L &= \Eb_{\tilde{\calL}} \Bigg[\bigg( -\sqrt{\rho} \int_{0}^{T} \Big\langle \vPhi(t,X,\vx ; \vTheta^{(k)}), \rd W(t) \Big\rangle\\
    & - \frac{1}{2} \rho \int_{0}^{T} \Big\langle \vPhi(t,X,\vx ; \vTheta^{(k)}), \vPhi(t,X,\vx ; \vTheta^{(k)}) \Big\rangle \rd t \bigg)  \frac{\rd  \calL^*}{\rd \calL} \frac{\rd  \calL}{\rd \tilde{\calL}} \rd \tilde{\calL} \Bigg]
\end{align*}
Evaluating, $\frac{\rd  \calL^*}{\rd \calL}\frac{\rd  \calL}{\rd \tilde{\calL}}$ separately, we have,
\begin{align*}
    \frac{\rd  \calL^*}{\rd \calL} \frac{\rd  \calL}{\rd \tilde{\calL}} &= \frac{\exp( - \rho J)}{\Eb_\calL \big[\exp( - \rho J) \big]} \, \exp(-\sqrt{\rho} \int_{0}^{T} \Big\langle \vPhi(t,X,\vx ; \vTheta^{(k)}), \rd W(t) \Big\rangle\\
    &- \frac{1}{2} \rho \int_{0}^{T} \Big\langle \vPhi(t,X,\vx ; \vTheta^{(k)}), \vPhi(t,X,\vx ; \vTheta^{(k)}) \Big\rangle \rd t)\\
    & = \frac{\exp( - \rho \tilde{J})}{\Eb_\calL \big[\exp( - \rho J) \big]}, 
\end{align*}
where $\tilde{J}$ is defined in \cref{eq:importance_sampled_cost}. Similarly, we can use importance sampling for the expectation in the denominator using \cref{eq:radon_i_nonlinear_feedback} as,
\begin{align*}
    \Eb_\calL \big[\exp( - \rho J) \big] &= \Eb_{\tilde{\calL}} \bigg[\frac{\rd  \calL}{\rd \tilde{\calL}} \exp( - \rho J) \bigg] = \Eb_{\tilde{\calL}} \big[\exp( - \rho \tilde{J}) \big]\\
    \therefore \frac{\rd  \calL^*}{\rd \calL} \frac{\rd  \calL}{\rd \tilde{\calL}} &= \frac{\exp( - \rho \tilde{J})}{\Eb_{\tilde{\calL}} \big[\exp( - \rho \tilde{J}) \big]}
\end{align*}
Putting all of this together, we get the required form of \cref{eq:loss_function} as,
\begin{align*}
     L &=  \Eb_{\tilde{\calL}}  \vast[\frac{\exp( - \rho \tilde{J})}{\Eb_{\tilde{\calL}} \big[ \exp( - \rho \tilde{J}) \big]} \bigg( -\sqrt{\rho} \int_{0}^{T} \Big\langle \vPhi(t,X,\vx ; \vTheta^{(k)}), \rd W(t) \Big\rangle\\
    & - \frac{1}{2} \rho \int_{0}^{T} \Big\langle \vPhi(t,X,\vx ; \vTheta^{(k)}), \vPhi(t,X,\vx ; \vTheta^{(k)}) \Big\rangle\rd t \bigg) \vast]    
\end{align*}

\section{Additional information on simulations}
Following are some details on each of our simulations which will help in reproducibility of the results. 

\subsection{Heat \ac{SPDE} distributed and boundary control}
The 2D Heat \ac{SPDE} with homogeneous Dirichlet boundary conditions given by: 
\begin{equation} \label{supeq:HeatSPDE}
\begin{split}
h_t(t, x, y) &= \epsilon h_{xx}(t,x,y) + \epsilon h_{yy}(t,x,y) + \sigma dW(t), \\
h(t,0,y) &= h(t,a,y) = h(t,x,0)= h(t,x,a)=0, \\ 
h(0,x,y) &\sim \mathcal{N}(h_0;0,\sigma_0), 
\end{split}
\end{equation}
where the parameter $\epsilon$ is the so called thermal diffusivity, which governs how quickly the initial temperature profile diffuses across the spatial domain. \eqref{supeq:HeatSPDE} considers the scenario of controlling a metallic plate to a desired temperature profile using 5 actuators distributed across the plate. The edges of the plate are always held at constant temperature of 0 degrees Celsius. The parameter $a$ is the length of the sides of the square plate, for which we use $a=0.25$ meters. 

The actuator dynamics are modelled by Gaussian-like exponential functions with the means co-located with the actuator locations at: $\vmu = \big[\mu_1, \mu_2, \mu_3, \mu_4, \mu_5 \big] =  \big[(0.2a,0.5a), (0.5a,0.2a), (0.5a,0.5a), \\(0.5a,0.8a), (0.8a,0.5a)\big]$ and the variance of the effect of each actuator on nearby field states given by $\sigma_l^2 = (0.1a)^2$, $\forall l = 1, \dots, 5$. The spatial domain is discretized by dividing the x and y domains into $J=32$ points each creating a grid of $32 \times 32$ spatial locations on the plate surface. The resulting $m_l(\vx)$ has the form:
\begin{align*}
    m_{l,j}\left(\left[\begin{array}{c} x \\ y\end{array} \right]\right) &= \exp \left\lbrace -\frac{1}{2}\left(\left[ \begin{array}{c} x \\ y
    \end{array} \right] - \left[ \begin{array}{c} \mu_{l,x} \\ \mu_{l,y} \end{array} \right] \right)^\top \left[ \begin{array}{cc}
         \sigma_l^2 & 0 \\
         0 & \sigma_l^2 
    \end{array} \right] \left(\left[ \begin{array}{c} x \\ y
    \end{array} \right] - \left[ \begin{array}{c} \mu_{l,x} \\ \mu_{l,y} \end{array} \right] \right) \right\rbrace, \\\quad \forall \, j &= 1, \dots, J, \;\; l=1, \dots, 5
\end{align*}

 For our simulations, we use a semi-implicit forward Euler discretization scheme for time and central difference for the $2^{nd}$ order spatial derivatives $h_{xx}$ and $h_{yy}$. We used the following parameter values, time discretization $\Delta t=0.02s$, simulation time horizon $T = 1.0s$ and thermal diffusivity $\epsilon=1.0$.
The cost function considered for the experiments was defined as follows: 
\begin{equation*}
    J := \sum_{t} \sum_{x} \sum_{y} \;\kappa \big(h_{\text{actual}}(t,x,y) - h_{\text{desired}} (t,x,y)\big)^2 \cdot \mathbbm{1}_{S}(x,y)
\end{equation*}
where $S := \cup_{i=1}^5 S_i$ and the indicator function $ \mathbbm{1}_{S}(x,y) $  is defined as follows:
\begin{align*}
\mathbbm{1}_S(x,y) :=
\begin{cases}
1,  \quad \text{if} ~~~(x,y) \in S  \\
0, \quad  \text{otherwise}
\end{cases}
\end{align*}
% \quad  \text{if}  ~~~0.48a  \leq x \leq 0.52 \;\text{and}\; 0.48a  \leq y \leq 0.52,   ~~~  \text{or} ~~~0.48  \leq x \leq 0.52   ~~~ \text{or} ~~~  0.78  \leq x \leq 0.82
where,\\ $S_1 = \{(x,y) \mid 0.48a  \leq x \leq 0.52a \;\text{and}\; 0.48a  \leq y \leq 0.52a\}$ is in the central region of the plate\\
$S_2 = \{(x,y) \mid 0.22a  \leq x \leq 0.18a \;\text{and}\; 0.48a  \leq y \leq 0.52a\}$ is the left-mid region of the plate \\
$S_3 = \{(x,y) \mid 0.82a  \leq x \leq 0.78a \;\text{and}\; 0.48a  \leq y \leq 0.52a\}$ is the right-mid region of the plate \\
$S_4 = \{(x,y) \mid 0.48a  \leq x \leq 0.52a \;\text{and}\; 0.18a  \leq y \leq 0.22a\}$ is in the top-central region of the plate \\
$S_5 = \{(x,y) \mid 0.48a  \leq x \leq 0.52a \;\text{and}\; 0.78a  \leq y \leq 0.82a\}$ is in the bottom-central region of the plate \\
\\
In addition  $ h_{\text{desired}} (t,x,y)= 0.5^{\circ}\,C $  for $(x,y) \in S_1$ and   $ h_{\text{desired}} (t,x,y)= 1.0^{\circ}\,C $  for $(x,y) \in S_2, S_3, S_4, S_5$ and the scaling parameter $\kappa=10^{-3}$.  

Since the domain is 2D, the inputs to the non-linear policy are image-like data and therefore the policy was chosen to be a \ac{CNN}. The description of the network 
architecture is as follows:

\begin{table}[H]
    \centering
    \begin{tabular}{| l | c |c |c| c| c|}
        \hline
        \textbf{Layer name} & \textbf{Kernel size} & \textbf{\# Filters (output size)} & \textbf{Stride} & \textbf{Padding type} & \textbf{Activation} \\
        \hline
        Input & - & 1 & - & - & -\\
        \hline
        Conv-1 & 4 & 5 & 2 & VALID & ReLU\\
        \hline
        Max-pool-1 & 2 & - & 2 & - & -\\
        \hline
        Conv-2 & 2 & 16 & 1 & SAME & ReLU\\
        \hline
        Max-pool-2 & 2 & - & 2 & - & -\\
        \hline
        Dense & - & 5 & - & - & Linear\\
        \hline
    \end{tabular} 
    \caption{Description of \ac{CNN} policy network for 2d Heat \ac{SPDE}.}
    \label{tab:CNN_network}
\end{table}

The network was trained using the ADAM optimizer for $1000$ iterations with $50$ trajectories ($1.0$ second long and $\Delta t=0.02$) sampled from the 2D Heat \ac{SPDE} model per iteration.

In the boundary control case, we make use of the 1D stochastic heat equation given as follows: 
\begin{align*}
h_t(t,x) &= \epsilon h_{xx}(t,x) +  \sigma dW(t)\\
h(0,x) &= h_0(x)
\end{align*}

For Dirichlet and Neumann boundary conditions we have $h(t,x) =\gamma(x)$, $\forall x\in\partial O$ and $h_x(t,x) =\gamma(x)$, $\forall x\in\partial O $, respectively. Regarding our 1-D boundary control example, we set $\epsilon=1 $, $\rho=10$, $h_x(t,0)=u_1(t) + \frac{1}{\sqrt{\rho}} dW(t)$ and $h_x(t,a)=u_2(t) + \frac{1}{\sqrt{\rho}} dW(t)$. In this case, $m_l(x)$ is simply given by an indicator function. Finally, the cost function used is the same as above with $S=\{x|0<x<a\}$ and $h_{desired}(t,x)= 3$. 

For 1D boundary control, the non-linear policy was chosen to be a \ac{FNN} with 2 hidden layers of 64 neurons each and ReLU activations. The network was trained using the ADAM optimizer for $1000$ iterations with $200$ trajectories (each $1.5$ seconds long and $\Delta t=0.01$ seconds) sampled from the Nagumo \ac{SPDE} model per iteration.

\subsection{1D Burgers \ac{SPDE} distributed control}
The 1D Burgers \ac{SPDE} with non-homogeneous Dirichlet boundary conditions is as follows:
\begin{equation} \label{supeq:BurgersSPDE}
\begin{split}
h_t(t, x) + h h_x(t, x) &= \epsilon h_{xx}(t,x) + \sigma dW(t)\\
h(t,0) &= h(t,a) = 1.0\\ 
h(0,x) &= 0, \; \forall x \in (0,a)
\end{split}
\end{equation}
where the parameter $\epsilon$ is the viscosity of the medium. \eqref{supeq:BurgersSPDE} considers a simple model of a 1D flow of a fluid in a medium with non-zero flow velocities at the two boundaries. The goal is to achieve and maintain a desired flow velocity profile at certain points along the spatial domain. As seen in the desired profile (Fig.3e) in the main paper, there are 3 areas along the spatial domain with desired flow velocity such that the flow has to be accelerated, then decelerated, and then accelerated again while trying to overcome the stochastic forces and the dynamics governed by the Burgers \ac{SPDE}. Similar to the experiments for the Heat \ac{SPDE}, we consider actuators behaving as Gaussian-like exponential functions with the means co-located with the actuator locations at: $\vmu = \big[0.2a, 0.3a, 0.5a, 0.7a, 0.8a\big]$ and the spatial effect (variance) of each actuator given by $\sigma_l^2 = (0.1a)^2$, $\forall \, l = 1, \dots, 5$. The parameter $a=1.0\;m$ is the length of the channel along which the fluid is flowing. 

This spatial domain was discretized using a grid of $64$ points. The numerical scheme used semi-implicit forward Euler discretization for time and central difference approximation for both the $1^{st}$ and $2^{nd}$ order derivatives in space. The $1^{st}$ order derivative terms in the advection term $u u_x$ were evaluated at the current time instant while the $2^{nd}$ order spatial derivatives in the diffusion term $u_{xx}$ were evaluated at the next time instant, hence the scheme is semi-implicit. Following are values of some other parameters used in our experiments: time discretization $\Delta t=0.01$, total simulation time = $1.0\,s$, and the scaling parameter $\kappa=100$. The cost function considered for the experiments was defined as follows: 
\begin{equation*}
    J := \sum_{t} \sum_{x} \;\kappa \big( h_{\text{actual}}(t,x) - h_{\text{desired}} (t,x)\big)^2 \cdot \mathbb{I}(x)
\end{equation*}
where the function  $ \mathbb{I}(x) $  is defined as follows 
$$
 \mathbb{I}(x) =
\begin{cases}
1, \quad  \text{if }  x \in 
\lbrace[0.18a,0.22a],[0.48a,0.52a],[0.78a,0.82a]\rbrace \\%~~~0.18a  \leq x \leq 0.22a,  ~~~  \text{or} ~~~0.48a  \leq x \leq 0.52a   ~~~ \text{or} ~~~  0.78a  \leq x \leq 0.82a \\
0, \quad  \text{otherwise}
\end{cases}
$$
In addition  $ h_{\text{desired}} (t,x)= 2.0 \;m/s$  for $x\;\in\lbrace[0.18a,0.22a],[0.78a,0.82a]\rbrace$ which is at the sides, and $ h_{\text{desired}} (t,x)= 1.0 \;m/s$  for $x\;\in [0.48a, 0.52a]$ which is in the central region.

The non-linear policy was chosen to be a \ac{FNN} with 2 hidden layers of 64 neurons each and ReLU activations. The network was trained using the ADAM optimizer for $1000$ iterations with $100$ trajectories (each $2.0$ seconds long and $\Delta t=0.01$ seconds) sampled from the Nagumo \ac{SPDE} model per iteration.

\subsection{1D Nagumo \ac{SPDE} distributed control (Suppression Task)}
\label{supsec:NagumoSPDE} The stochastic Nagumo equation with Neumann boundary conditions is as follows: 
\begin{align*}
h_t(t,x) &= \epsilon h_{xx}(t,x) + h(t,x)\big(1-h(t,x)\big)\big(h(t,x)-\alpha\big) + \sigma dW(t)\\ h_x(t,0) &= h_x(t,a) = 0\\ 
h(0,x) &= \bigg(1+\exp\Big(-\frac{2-x}{\sqrt[]{2}}\Big)\bigg)^{-1}
\end{align*}
The parameter $\alpha$ determines the speed of a wave traveling down the length of the axon and $\epsilon$ the rate of diffusion. By simulating the deterministic Nagumo equation with $a=5.0,\,\epsilon=1.0$ and $\alpha=-0.5$, we observed that after about 3.5 seconds, the wave completely propagates to the end of the axon. We consider actuators behaving as Gaussian-like exponential functions with actuator centers (mean values) at $\vmu = \big[0.7a, 0.8a, 0.9a\big]$ and the spatial effect (variance) of each actuator given by $\sigma_l^2 = (0.1a)^2$, for $\, l = 1, 2, 3$. The spatial domain was discretized using a grid of 64 points. The cost function for this experiment was defined as follows: 
\begin{equation*}
    J = \sum_{t} \sum_{x} \;\kappa \big(h_{\text{actual}}(t,x) (t,x)\big)^2 \cdot \mathbb{I}(x)
\end{equation*}
where $\kappa$ was chosen as $10^{-3}$, and the function  $ \mathbb{I}(x) $  is defined as follows 
$$
 \mathbb{I}(x) =
\begin{cases}
1, \quad  \text{if}  ~~~ x\in [0.7a, 0.99a] \\% \leq x \leq 0.99a \\ %\geq 0.7a  \;\text{and} \;x \leq 0.99a\\
0, \quad  \text{otherwise}
\end{cases}.
$$
The non-linear policy was chosen to be a \ac{FNN} with 2 hidden layers of 64 neurons each and ReLU activations. The network was trained using the ADAM optimizer for $1000$ iterations with $50$ trajectories (each $3.5$ seconds long and $\Delta t=0.01$ seconds) sampled from the Nagumo \ac{SPDE} model per iteration.

\end{document}